\documentclass[12pt, a4paper,english]{article}
\usepackage[T1]{fontenc}
\usepackage[utf8]{inputenc}
\usepackage{babel}
\usepackage{url}

\usepackage{amsfonts, amsthm, amsmath, amssymb}
\usepackage{textcomp, vmargin, microtype}

\usepackage[]{hyperref}

\newtheorem{teo}{Theorem}[section]
\newtheorem{theorem}[teo]{Theorem}

\theoremstyle{definition}

\numberwithin{equation}{section}

\usepackage{tikz,xcolor,hyperref}
\definecolor{lime}{HTML}{A6CE39}
\DeclareRobustCommand{\orcidicon}{%
	\begin{tikzpicture}
	\draw[lime, fill=lime] (0,0) 
	circle [radius=0.16] 
	node[white] {{\fontfamily{qag}\selectfont \tiny ID}};
	\draw[white, fill=white] (-0.0625,0.095) 
	circle [radius=0.007];
	\end{tikzpicture}
	\hspace{-2mm}
}
\foreach \x in {A, ..., Z}{%
	\expandafter\xdef\csname orcid\x\endcsname{\noexpand\href{https://orcid.org/\csname orcidauthor\x\endcsname}{\noexpand\orcidicon}}
}

\title{The geodesics for Poincar\'e's half-plane: 
a nonstandard derivation}

\author{
Gianluca Gorni \orcidB{}\\
Universit\`a di Udine,
Dipartimento di\\ Scienze Matematiche, Informatiche e Fisiche\\
via delle Scienze~208, 33100 Udine, Italy\\
\tt{gianluca.gorni@uniud.it}
\and
Gaetano  Zampieri \orcidC{}\\
Universit\`a di Verona,
Dipartimento di Informatica\\
strada Le Grazie 15, 37134 Verona, Italy\\
\tt{gaetano.zampieri@univr.it}}

\date{}
\begin{document}
\maketitle

\begin{abstract}
Constants of motion are usually derived from groups of symmetry transformation of the system. Here we show that useful properties of the system can be deduced from a family of Noether-like transformations that are not inspired by any symmetry whatsoever. The system here is the Lagrangian interpretation of Poincar\'e's half plane, and the property is the shape of the geodesics.
\end{abstract}

{\textbf{Keywords:} }

\begin{center}
Dedicated to our thesis advisors Giuseppe Da Prato and Aldo Bressan
\end{center}

\begin{center}
This is the original manuscript, version of October 1, 2021
\end{center}
\section{Introduction}

Lobachev\-sky's hyperbolic geometry was given models by Beltrami, and then Klein and Poincar\'e. Here we are interested in the model that is called Poincar\'e's half plane: the set $D=\{(x,y)\in\mathbb{R}^2: y>0\}$ where we call ``straight line'' any half-circle with center on the $x$ axis and any vertical half-line. As in Fig.~\ref{poincareHalfPlane}, given a ``straight line''~$r$ (the thicker half-circle) and a point~$P$, there are infinitely many ``straight lines'' that cointain~$P$ and do not intersect~$r$.

Poincar\'e's half plane can be formulated in the language of differential geometry, because the ``straight lines'' coincide with the geodesics with respect to the metric $(\begin {smallmatrix} 1&0\\ 0&1\end{smallmatrix}) /y^2$. Then there is the point of view of variational mechanics, on which we will concentrate: the ``straight lines'' are the variational stationary points, or natural motions, for the action defined by the following Lagrangian function:
\begin{equation}\label{LPoincare}
  L(t,q,\dot q)=\frac{\dot q_1^2+\dot q_2^2}{2q_2^2}\,,
  \qquad t\in\mathbb{R}, \quad q\in D,\quad 
  \dot q\in\mathbb{R}^2,
\end{equation}
where we switch to the notation that is familiar in mechanics: $(x,y)=q=(q_1,q_2)$. We were attracted to this Lagrangian formulation because we are in search for applications of our theory of ``nonlocal constants of motion''. In~so doing we happened to fine an unusual derivation of the well-known fact that the nonconstant motions for~$L$ are indeed half-circles and vertical half-lines. The calculations are simple enough to be presented to a general audience, while the specialist may find it surprising that useful constants of motion can be found without appealing to Noether symmetries.

\section{Nonlocal constants of motion}

The natural motions for a Lagrangian~$L$ are the solutions to the Euler-Lagrange equations 
\begin{equation}\label{EL}
  \frac{d}{dt}\frac{\partial L}{\partial\dot q}
  \bigl(t,q(t),\dot q(t)\bigr)-
  \frac{\partial L}{\partial q}\bigl(t,q(t),\dot q(t)\bigr)=0,
\end{equation}
which, for our Lagrangian \eqref{LPoincare} in Poincar\'e's half-plane, become
\begin{equation}\label{LEPoincare}
  \ddot q_1-\frac{2}{q_2}\,\dot q_1\,\dot q_2=0,\qquad
  \ddot q_2+\frac{1}{q_2}\bigl(\dot q_1^2-\dot q_2^2\bigr)=0.
\end{equation}
This Lagrangian system has the following functions $E,p$ as immediate first integrals
\begin{equation}\label{PoincareFirstIntegrals}
  \frac{\dot q_1^2+\dot q_2^2}{2q_2^2}=E,\qquad
  \frac{\dot q_1}{q_2^2}=p\,,
\end{equation}
The first is \emph{energy}, which is a property of any 
time-independent Lagrangian $L(t,q,\dot q)=\mathcal{L}(q,\dot q)$ through the formula
\begin{equation}
  E=\frac{\partial {\cal L}}{\partial\dot q}
  \bigl(q,\dot q\bigr)\cdot \dot q-{\cal L}\bigl(q,\dot q\bigr)
\end{equation}
(where the central dot is the scalar product in $\mathbb{R}^n$). The function~$p=\dot q_1/q_2^2$ (the first component of \emph{momentum}), is constant because of the lack of $q_1$ in the Lagrangian. Indeed, from the Euler-Lagrange equations~\eqref{EL}
\begin{equation}
  \frac{\partial L}{\partial q_1}
  \bigl(t,q,\dot q\bigr)\equiv 0\quad \Longrightarrow \quad
  \frac{\partial L}{\partial\dot q_1}
  \bigl(t,q(t),\dot q(t)\bigr)=\text{constant}.
\end{equation}
Clearly, $E>0$ for nonconstant geodesics.

In the spirit of Noether's theorem, energy conservation is a consequence of the invariance of the Lagrangian by \emph{time-translations}, and the other first integral comes from invariance under \emph{$q_1$-translations}. In a previous paper \cite[subsection 7.2.2]{GZsao} we also used one more first integral, which comes from invariance under \emph{dilations}, but here we pointedly avoid~it.

What we do is to consider \emph{$q_2$-translations}. This may seem ludicrous at first sight, because $q_2$-translations do \emph{not} leave $L$ invariant. Nevertheless, it makes perfect sense within our framework that generates nonlocal constants of motion, which can be stated simply enough:

\begin{theorem}\label{nonlocaltheorem2nd} Let $L(t,q,\dot q)$ be a smooth scalar valued Lagrangian function,  $t\in \mathbb{R}$, $q,\dot q\in\mathbb{R}^n$. Let $q(t)$ be a solution to the Euler-Lagrange equation and let $q_\lambda(t)$, $\lambda\in \mathbb{R}$, be a smooth family of perturbed  motions, such that $q_0(t)\equiv q(t)$.
Then the following function of $t$ is constant
\begin{equation}\label{constantofmotionsecondorder}
 \frac{\partial L}{\partial \dot q}
  \bigl(t,q(t),\dot q(t)\bigr)\cdot
  \frac{\partial q_\lambda}{\partial\lambda} (t)
  \Big|_{\lambda=0} -
  \int_{t_0}^t  \frac{\partial}{\partial\lambda}
  L\bigl(s,q_\lambda(s),\dot q_\lambda(s)\bigr)
  \big|_{\lambda=0}ds\,.
\end{equation}
\end{theorem}

\noindent
The proof is straightforward: we just take the derivative of the function in~\eqref{constantofmotionsecondorder} and use the Euler-Lagrange equation and reverse the order of a double derivative.

This expression \eqref{constantofmotionsecondorder} will be called \emph{the constant of motion associated to the family} $q_\lambda(t)$. For a random family, we may expect the constant of motion to be trivial or inconsequential. In general it is \emph{nonlocal}, which means that its value at a time $t$ depends not only on the current state $(t,q(t),\dot q(t))$ at time~$t$, but also on the whole history between an arbitrary $t_0$ and~$t$.

In the original spirit of Noether's theorem we would concentrate the attention to families $q_\lambda(t)$ which make the integrand in \eqref{constantofmotionsecondorder} vanish whenever $L$ enjoys an invariance property. For instance, the $q_1$-translation family $q_{\lambda}(t) =(q_1(t)+\lambda,q_2(t))$, if plugged into formula~\eqref{constantofmotionsecondorder} yields the first integral~$p$ in \eqref{PoincareFirstIntegrals} for the geodesics of Poincar\'e's half-plane.

As we are going to see in the next Section~\ref{separation}, with $q_2$-translations the integral of formula~\eqref{constantofmotionsecondorder} does not disappear, but, in spite of that, we will put the nonlocal constant of motion to good use, as it allows for a nonstandard separation of variables.

Theorem~\ref{nonlocaltheorem2nd} is part of a line of research that started in 2014, the interested reader may find applications and references in~\cite{GZhydraulicdissipative}, in particular Section~6.1 deals with a similar nonstandard separation of variables in dimension 6 for the conservative Maxwell-Bloch equations. 
\section{Integration of the geodesics}\label{separation}

Consider the Lagrangian~\eqref{LPoincare}. 
As announced, let us apply Theorem~\ref{nonlocaltheorem2nd} with the \emph{$q_2$-translation} family $q_{\lambda}(t) =\bigl( q_1(t), q_2(t) +\lambda\bigr)$. Then
\begin{equation*}
  \frac{\partial q_\lambda}{\partial\lambda}(t)
  \big|_{\lambda=0}=\bigl(0,1\bigr),\quad  
   \frac{\partial }{\partial\lambda}
  L\bigl(t,q_\lambda(t),
  \dot q_\lambda(t)\bigr)\big|_{\lambda=0}=
  -\frac{\dot q_1(t)^2+\dot q_2(t)^2}{q_2(t)^3}.
\end{equation*}
and the nonlocal constant of motion becomes
\begin{equation}\label{nonlocalConstantForHalfPlane}
 \frac{\dot q_2(t)}{q_2(t)^2} +
  \int_{t_0}^t
  \frac{\dot q_1(s)^2+\dot q_2(s)^2}{q_2(s)^3} ds=
  -\frac{d}{dt}\frac{1}{q_2(t)} +
  2E\int_{t_0}^t \frac{1}{q_2(s)} ds
\end{equation}
where we have used the energy conservation, i.e., the fact that the~$E$ in~\eqref{PoincareFirstIntegrals} is constant along the motions. Since \eqref{nonlocalConstantForHalfPlane} is constant along any solution to the Euler-Lagrange equation, its time derivative vanishes:
\begin{equation}\label{equationIn1/q2}
  -\frac{d^2}{dt^2}\frac{1}{q_2(t)}+
  2E\frac{1}{q_2(t)}=0.
\end{equation}
This separate linear equation in the new dependent variable $1/q_2$ is integrated at once. If~$E=0$ we already know that $(q_1,q_2)$ is constant. If~$E>0$ the generic solution of~\eqref{equationIn1/q2} is
\begin{equation}\label{formulaForQ2}
  q_2(t)=\left(c_1 e^{t\sqrt{2E}}+
  c_2  e^{-t\sqrt{2E}}\right)^{-1}.
\end{equation}
Since $q_2(0)>0$ it must be that $c_1+c_2>0$. Then eliminate $\dot q_1$ from the conservation laws in~\eqref{PoincareFirstIntegrals}, apply formula \eqref{formulaForQ2} and set $t=0$: this way we obtain a necessary condition for $c_1,c_2$ to give a geodesic:
\begin{equation}\label{conditionForGeodesic}
  8c_1c_2E=p^2.
\end{equation}

If $p=0$ the geodesic is contained in a vertical half-line because the conservation law $\dot q_1=pq_2^2$ implies that $q_1(t)$ is constant; also, \eqref{conditionForGeodesic} gives that either $c_1=0$ or $c_2=0$. The geodesic then becomes fully explicit as $q_2(t)= e^{\pm t\sqrt{2E}}/c$ with $c=1/q_2(0)>0$. In either case $t$ ranges in the whole of~$\mathbb{R}$. 

In the nondegenerate case $p\ne0$, both $c_1>0$ and $c_2>0$, $q_2(t)>0$ is defined for all $t\in\mathbb{R}$, and
we can deduce $q_1$ using $\dot q_1=pq_2^2$ up to an integration constant~$c_3$:
\begin{equation}\label{formulaForQ1}
  q_1(t)= c_3-
  \frac{p}{2c_1\sqrt{2E}
  \bigl(c_2+c_1e^{2t\sqrt{2E}}\bigr)},
\end{equation}
also defined for all $t\in\mathbb{R}$. Let us write $x=q_1(t)$, $y=q_2(t)$ and eliminate the exponentials from formulas \eqref{formulaForQ2} and~\eqref{formulaForQ1}: we obtain the relation
\begin{equation}
  \frac{8c_1c_2E}{p^2}\Biggl(x-c_3
  +\frac{p}{2c_1c_2\sqrt{8E}}\Biggr)^2
  +y^2=\frac{1}{4c_1c_2},
\end{equation}
which is the equation of an ellipse. If we finally impose the necessary condition~\eqref{conditionForGeodesic}, the ellipse becomes a circle
\begin{equation}
  (x-c)^2+y^2=\frac{2E}{p^2},
\end{equation}
with center on the $x$~axis and radius $\sqrt{2E}/\lvert p\rvert$.

We have proved that the geodesics of the Poincar\'e half-plane are either half-lines or half-circles using the two conservation laws for~$E$ and~$p$ plus the nonlocal constant of motion~\eqref{nonlocalConstantForHalfPlane}.

\begin{figure}
  \begin{center}
\includegraphics[width=.8\textwidth]{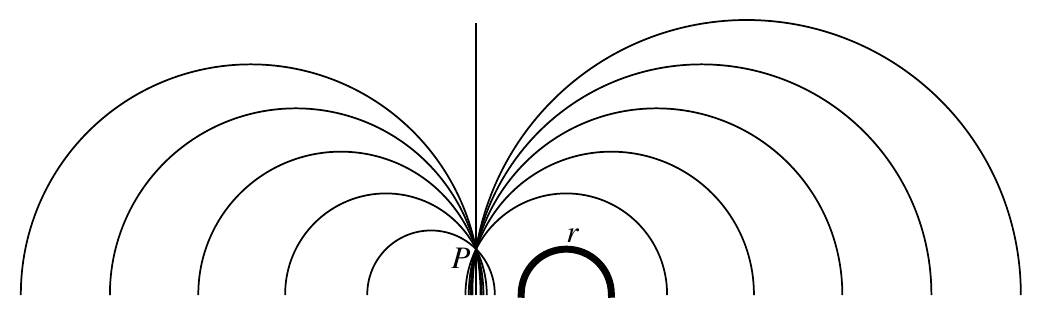}
\end{center}
\caption{Geodesics in Poincar\'e's half-plane.}
\label{poincareHalfPlane}
\end{figure}

On Poincar\'e's half-plane every couple of points is connected by a geodesic. Figure~\ref{poincareHalfPlane} shows that given a geodesic~$r$ and a point~$P$ outside, there exist infinitely many geodesics through~$P$ that do not cross~$r$.

\goodbreak
\section*{Acknowledgments}

Paper written under the auspices of INdAM (Istituto Nazionale di Alta Matematica).



\begin{thebibliography}{10}

\bibitem{Beltrami}
Beltrami, E. (1868).
\newblock Teoria fondamentale degli spazi di curvatura constante.
\newblock {\em Annali di Matematica Pura ed Applicata}, ser II 2, 232--255.
\url{https://gallica.bnf.fr/ark:/12148/bpt6k99432q}

\bibitem{GZsao}
Gorni, G.,  Zampieri, G. (2011).
\newblock Variational aspects of analytical mechanics.  
\newblock {\em S\~ao  Paulo J. Math. Sci.} 5: 203--231.
DOI: 10.11606/issn.2316-9028.v5i2p249-279


\bibitem{GZhydraulicdissipative}
Gorni  G., Zampieri G. (2020).
\newblock Lagrangian dynamics by nonlocal constants of motion.
\newblock {\em Discrete Contin. Dyn. Syst. Ser. S} 13(10): 2751--2759. DOI: 10.3934/dcdss.2020216


\bibitem{Poincare}
Poincar\'e, H. (1882).
\newblock Th\'eorie des Groupes Fuchsiens. 
\newblock{\em Acta Mathematica}.
\url{https://archive.org/details/thoriedesgroup00poin/page/n67}

\end{thebibliography}
\end{document}